\documentclass[10pt]{amsart}
\usepackage{amsfonts,amssymb,amscd}
\usepackage{pstricks}
\usepackage{pst-plot}
\usepackage{mathrsfs}

\let\mathcal\mathscr

\usepackage[all,ps,cmtip,color]{xy}

\newtheorem{theo}{Theorem}

\newtheorem{prop}{Proposition}

\begin{document}

\def\om{\omega}
\def\ot{\otimes}
\def\we{\wedge}
\def\wec{\wedge\cdots\wedge}
\def\op{\oplus}
\def\ra{\rightarrow}
\def\lra{\longrightarrow}
\def\fso{\mathfrak so}\def\fg{\mathfrak g}
\def\cO{\mathcal{O}}
\def\cS{\mathcal{S}}
\def\fsl{\mathfrak sl}\def\fsp{\mathfrak sp}
\def\fb{\mathfrak b}
\def\fh{\mathfrak h}
\def\fn{\mathfrak n}
\def\PP{\mathbb P}\def\QQ{\mathbb Q}\def\ZZ{\mathbb Z}\def\CC{\mathbb C}
\def\RR{\mathbb R}\def\HH{\mathbb H}\def\OO{\mathbb O}\def\AA{\mathbb A}
\def\AP2{{\mathbb{AP}^2}}
\def\RP2{{\mathbb{RP}^2}}\def\RPC{{\mathbb{RP}_\CC^2}}
\def\CP2{{\mathbb{CP}^2}}
\def\HP2{{\mathbb{HP}^2}}
\def\APn{{\mathbb{AP}^n}}\def\APnC{{\mathbb{AP}_\CC^n}}
\def\RPn{{\mathbb{RP}^n}}
\def\CPn{{\mathbb{CP}^n}}
\def\HPn{{\mathbb{HP}^n}}
\def\OP2{{\mathbb{OP}^2}}
\def\APC{{\mathbb{AP}^2_\mathbb{C}}}
\def\JA{J_3(\mathbb{A})}
\def\JC{J_3(\mathbb{A}_{\mathbb{C}})}\def\JH{J_3(\mathbb{H}_{\mathbb{C}})}
\def\JnC{J_{n+1}(\mathbb{A}_{\mathbb{C}})}
\def\smc{\cdots }
\def\Gr{{\mathcal Gr}}
\def\Gm{\mathbf{G_m}}

\title{Triangulations and Severi varieties}

\author{F. Chapoton, L. Manivel}


\begin{abstract}
  We consider the problem of constructing triangulations of projective
  planes over Hurwitz algebras with minimal numbers of vertices. We
  observe that the numbers of faces of each dimension must be equal to
  the dimensions of certain representations of the automorphism groups
  of the corresponding Severi varieties.  We construct a complex
  involving these representations, which should be considered as a
  geometric version of the (putative) triangulations.
\end{abstract}

\maketitle

\section{Introduction}

Compare the following two statements, one from complex projective geometry, 
the other one from combinatorial topology. 

\begin{theo}[Zak, 1982]
 Let $X^d\subset \PP^{N-1}$ be a smooth irreducible complex projective variety of dimension $d$. 
\begin{enumerate}
 \item If $N<3\frac{d}{2}+3$, then the secant variety of $X$ fills out the ambient
space, $Sec(X)=\PP^{N-1}$.
 \item If $N=3\frac{d}{2}+3$, then either $Sec(X)=\PP^{N-1}$, or $d=2,4,8,16$.
\end{enumerate}
\end{theo}

Recall that the secant variety $Sec(X)$ is obtained by taking the union of the 
lines joining any two points of $X$, and passing to the Zariski closure.
  
The only exceptions to the second statement are the {\it Severi varieties}, 
the complexifications $\APC$ of the projective planes over  
$\AA=\RR, \CC, \HH, \OO$, the four normed algebras (see e.g. \cite{baez}). 

\begin{theo}[Brehm-K\"uhnel, 1987]
Let $X^d$ be a combinatorial manifold of dimension $d$, having $N$ vertices. 
\begin{enumerate}
 \item If $N<3\frac{d}{2}+3$, then $X$ is topologically a sphere. 
 \item If $N=3\frac{d}{2}+3$, then either $X$ is a sphere, or $d=2,4,8$, or $16$.
\end{enumerate}
\end{theo}

Possible exceptions to the second statement are the {\it real Severi varieties}, 
the projective planes $\AP2 $ over  $\AA=\RR, \CC, \HH, \OO$.

More precisely, there is a classical triangulation of the real projective 
plane $\RP2$ with $6$-vertices, described in the picture below where opposite sides 
of the big triangle must be identified. 

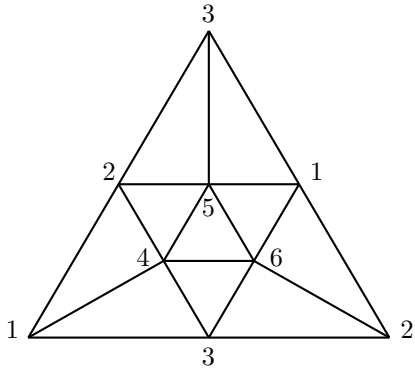
\begin{center}
\psset{unit=3mm}
\psset{xunit=3mm}
\psset{yunit=5.1mm}
\begin{pspicture*}(20,9)(-3,-1)
\psline(0,0)(16,0)\psline(0,0)(8,8)\psline(8,8)(16,0)
\psline(0,0)(6,2)\psline(8,4)(8,8)\psline(10,2)(16,0)
\psline(4,4)(12,4)\psline(4,4)(8,0)\psline(8,0)(12,4)
\psline(6,2)(10,2)\psline(10,2)(8,4)\psline(8,4)(6,2)
\put(-1,0){$1$}\put(3.3,7){$2$}\put(7.7,14){$3$}
\put(12.5,7){$1$}\put(16.5,0){$2$}\put(7.7,-1.2){$3$}
\put(4.8,3.2){$4$}\put(7.7,5.4){$5$}\put(10.7,3.2){$6$}
\end{pspicture*}
\end{center}

There is a unique triangulation of the 
complex projective plane $\CP2$ with only $9$-vertices \cite{bak}. Over the
quaternions the situation is not completely clear: it was shown in \cite{bk2}
that there exists three different combinatorial triangulations with 15 vertices 
of an eight-dimensional manifold which is ``like the quaternionic projective plane'',
but the authors could not decide whether this topological manifold was indeed $\HP2$
or a fake quaternionic plane. Finally, there is no candidate for a combinatorial 
triangulation of $\OP2$ with 27 vertices. We will see that the number of maximal
faces of such a triangulation should be $100386$!

Recall that the homology of $\AP2$ is 
$$H_i(\AP2, k) = \Bigg\{ \begin{array}{ll} 
                                  k & \mathrm{if}\; i=0,a,2a, \\
                                   0 & \mathrm{otherwise},
                                 \end{array} $$
if $k$ stands for $\ZZ$ when $a\ge 2$, and for $\ZZ_2$ when $a=1$. 
A more uniform statement is that each $\AP2$ has a Morse function with only 
three critical points. The precise statement of the second assertion of 
the theorem of Brehm and K\"uhnel is that $X$, if not a sphere, must 
admit such a Morse function. Manifolds with this property, which are 
{\it like projective planes}, were studied systematically
in \cite{ek}. Among other topological restrictions, 
the fact that the dimension of such a manifold 
must be $2,4,8$ or $16$ is established there. 

The goal of our paper is to explore the relationships between the two
statements above. Our main observation will be that the numbers of
faces of each dimension in a (putative) triangulation of a projective
plane $\AP2$, must be equal to the dimensions of certain linear
representations of the automorphism groups of the corresponding Severi
varieties $\APC$. Moreover, we will construct complexes involving
these representations, which we conjecture to be closely related with
the face complexes of the triangulations. Over the complex and
quaternionic numbers we will reconsider the results of \cite{bk1} and
\cite{bk2} and check that they correctly fit with our perspective. We
then present an observation concerning the links of vertices in the
triangulations. In a final section we elaborate on possible extensions
to projective spaces of higher dimensions.

\section{The Severi varieties}

We briefly recall the main geometric properties of the Severi varieties 
$\APC$ (see e.g. \cite{baez} and references therein). 
In all the sequel we denote by $a=1,2,4,8$ the dimension of 
the Hurwitz algebra 
$\AA=\RR, \CC, \HH, \OO$ as a real vector space. Recall that each 
$\AA$ has a natural 
involution generalizing the usual complex conjugation; it can be 
defined as the orthogonal symmetry with respect to the unity. 
Then consider the 
space $\JA$ of Hermitian $3\times 3$ matrices with coefficients
in $\AA$ . 
This is a real vector space of dimension $3a+3$, endowed with a 
structure of Jordan algebra defined by symmetrization of the 
ordinary  matrix multiplication. 

The automorphism group of the Jordan algebra $\JA$ will be denoted $SO(3,\AA)$. 
It preserves the cubic form defined by the determinant, which 
exists even over the octonions. The group of invertible linear 
transformations of the vector space $\JA$ preserving this determinant 
will be denoted  $SL(3,\AA)$. 

Let $\JC$ denote the complexification of $\JA$. The Severi variety 
$\APC\subset \PP \JC$ can be defined as the cone over the set of 
rank one matrices, where having rank one is defined by the property 
that all the derivatives of the determinant vanish (these derivatives 
are the analogues of $2\times 2$ minors). Geometrically, 
this means that the Severi variety is the singular locus of the 
determinantal hypersurface. We will only mention a few of its many 
remarkable properties:
\begin{enumerate}
 \item $\APC$ is smooth of dimension $2a$. 
 \item $\APC$ is homogeneous under the action of the complexified 
group $SL(3,\AA_\CC)$. 
 \item The secant variety of $\APC$ is the determinantal 
hypersurface. (The sum of two rank one matrices has rank at most two!)
In particular $\APC$ is {\it secant defective}, in the sense that
the secant variety is smaller that expected. 
\end{enumerate}

Polarizing the determinant, one obtains a quadratic map $c: \JC\rightarrow \JC^\vee$ 
which we call the {\it comatrix map}. Let $$W(\AA)=\CC\oplus \JC\oplus \JC^\vee\oplus\CC.$$
One can define the projective variety $LG(3,\AA_\CC)$ as the image of the rational map from $\JC$
to $\PP W(\AA)$ mapping $x\in\JC$ to the line generated by $1+x+c(x)+\det(x)$.
Denote by $p$ the point of $LG(3,\AA_\CC)$ defined by $x=0$.
The following properties do hold:
\begin{enumerate}
 \item $LG(3,\AA_\CC)$ is smooth of dimension $3a+3$. 
 \item $LG(3,\AA_\CC)$ is homogeneous under the action of a simple 
Lie group $Sp(6,\AA_\CC)$. 
 \item The lines through $p$ contained in $LG(3,\AA_\CC)$ generate a cone 
over $\APC$. 
\end{enumerate}

The groups we met have the following types:
$$\begin{array}{rcccc}
 \AA & \RR & \CC & \HH & \OO \\
SO(3,\AA_\CC) & A_1 & A_2 & C_3 & F_4 \\
SL(3,\AA_\CC) & A_2 & A_2\times A_2 & A_5 & E_6 \\
Sp(6,\AA_\CC) & C_3 & A_5 & D_6 & E_7
  \end{array}$$
This table is a chunk of the famous Tits-Freudenthal magic square. 
(For more on this see the survey paper \cite{LMsur} and references therein.)
\medskip

\section{Faces and representations}

\subsection{Special properties of minimal triangulations}
The triangulations of $\RP2$, $\CP2$ and (supposedly) $\HP2$
with minimal numbers of vertices 
have very peculiar properties \cite{bk1, bk2}, among which: 
\begin{enumerate}
\item (Tightness) Each face of dimension $a$ or less is part of the triangulation.
\item (Duality) A face of dimension $a+j$ is part of the triangulation if and 
only if the complementary face of dimension $2a-j+1$ is not.
\item (Secant defectivity) Any two maximal simplices intersect along a simplex
of dimension at least $a-1$. 
\end{enumerate}
 
Note that tightness is a consequence of duality, since there is no face
of dimension $2a+1$. Moreover, it was noticed by Marin that the duality
property is imposed by the algebra structure of the $\ZZ_2$-valued 
cohomology (see \cite{am}, and \cite[Proposition 2]{bk2}).

Since the number of vertices is $3a+3$, the dimension of the
intersection of two simplices of dimension $2a$ is at least $a-2$, and
should in general be equal to $a-2$.  In our triangulations any two
maximal simplices meet in dimension $a-1$, and this means that their
linear span has dimension $3a+1$ rather than the expected $3a+2$. This
is why this property should be understood as the combinatorial version
of the secant defectivity property of the Severi varieties.

\subsection{Numbers of faces}
The numbers of faces $f_k$ of each dimension $k$ in a triangulation of a 
smooth manifold are not independent. For example a codimension one face 
has to belong to exactly two codimension zero faces, hence the relation 
$$(d+1)f_d=2f_{d-1}$$ if $d$ denotes the dimension. More generally, 
the Dehn-Sommerville equations for simplicial complexes, or rather their 
extension by Klee (see \cite[Theorem 5.1]{ns}) to a context including 
combinatorial manifolds, imply that the numbers of faces of dimension 
smaller than half of $d$ determine the remaining numbers of faces. 
The precise statement is the following. Let $f_{-1}=1$, and consider
the $h$-vector, which is the sequence $h_0, \ldots , h_{d+1}$ defined 
by the identity
$$\sum_{i=0}^{d+1}h_ix^{d+1-i}=\sum_{j=0}^{d+1}f_{j-1}(x-1)^{d+1-j}.$$
For a simplicial complex, the classical Dehn-Sommerville equations
assert that the $h$-vector is symmetric, that is $h_i=h_{d+1-i}$. 
More generally, for a triangulation of a 
smooth manifold $X$, the $h$-vector is such that  
$$h_{d+1-i}-h_i=(-1)^i\binom{d+1}{i}\Big( \chi_{top}(X)-\chi_{top}(S^d)\Big) ,$$
where $\chi_{top}(X)$ denotes the topological Euler characteristic.
 
Let us denote by $f_k^a$ the number of faces of dimension $k$ in a tight 
triangulation of $\AP2$ with $3a+3$ vertices. The tightness property means that 
$$f_k^a=\binom{3a+3}{k+1} \qquad \mathrm{for}\quad 0\le k\le a.$$
The generalized Dehn-Sommerville equations show that all the 
numbers $f_k^a$ are then completely determined. 
These numbers are given in the following table:

$$\begin{array}{lllll}
& \RP2\; & \CP2\;\; & \HP2\;\;\; & \OP2 \\
\# \mathrm{vertices} & 6 & 9  & 15 &  27 \\
& 15 & 36 & 105 & 351 \\
\# \mathrm{2-dim} & 10 & 84 & 455 & 2925 \\
&   & 90 & 1365 & 17550 \\ 
\# \mathrm{4-dim}&   & 36 & 3003 & 80730 \\
&   &    & 4515 & 296010 \\
&   &    & 4230 & 888030 \\
&   &    & 2205 & 2220075 \\
\# \mathrm{8-dim}&   &    & 490  & 4686825 \\
&  &&& 8335899\\
&  &&& 12184614\\
&  &&& 14074164\\
&  &&& 12301200\\
&  &&& 7757100\\
&  &&& 3309696\\
&  &&& 853281 \\
\# \mathrm{16-dim} &  &&& 100386
\end{array}
$$

\bigskip\noindent {\bf Main observation.} {\it For each $k$, the number $f_k^a$ of 
$k$-dimensional faces of a minimal triangulation of $\AP2$ is the dimension of a 
representation of $SL(3,\AA_\CC)$.} 

\subsection{A connection with Severi varieties}
One can be much more precise. We will describe a recipe which allows to 
understand a priori which representation of $SL(3,\AA_\CC)$ has dimension $f_k^a$. 
Note that it will not be an irreducible representation in general, but 
it will have very few irreducible components, and never more than three. 

\medskip
 Consider the  variety $LG(3,\AA_\CC)=G/P$,
where $G=Sp(6,\AA_\CC)$ and $P$ is the stabilizer of the base point $p$. As any 
rational homogeneous variety does, $LG(3,\AA_\CC)$ has a cellular decomposition 
defined by the Schubert cells. If $B\subset P\subset G$ is a Borel subgroup,
recall that the Schubert cells can be defined as the $B$-orbits inside $LG(3,\AA_\CC)$. 
Their closures are called Schubert varieties and usually denoted $X_u$, 
where $u$ belongs to some index set $W_P$ defined in terms of the combinatorics
of the root system of $G$. It is clear from the 
definition that the boundary of any Schubert variety is a finite union
of smaller Schubert varieties. This allows to define an oriented graph, 
which we call the {\it Hasse diagram}. The vertices of this graph are in bijection 
with the Schubert varieties, that is with $W_P$. Moreover, there is an arrow
$u\rightarrow v$ if $X_u$ is an irreducible component of the boundary 
of $X_v$ (or equivalently, $X_u$ is a codimension one subvariety of $X_v$).  
The Hasse diagram is obviously ranked by the dimension of the
Schubert varieties. Moreover, Poincar\'e duality implies that it is 
symmetric with respect to the middle dimension. We denote the operation 
of Poincar\'e duality on the Hasse diagram by $\pi$.

\begin{prop}
One has the following properties:
\begin{enumerate}
 \item The cone over $\APC$ is a Schubert variety of $LG(3,\AA_\CC)$. 
 \item The Hasse diagram of $\APC$ embeds in the Hasse diagram of 
$LG(3,\AA_\CC)$, as an interval $I$. 
 \item The Hasse diagram of $LG(3,\AA_\CC)$ is the disjoint union 
of the interval $I$, the Poincar\'e dual $\pi(I)$, and the two 
extremities given by the fundamental class and the punctual class.  
\end{enumerate}
 \end{prop}

\proof The first claim is clear since the cone over $\APC$ is the union 
of the lines in $LG(3,\AA_\CC)$ passing through $p$. In particular it
is stabilized by $P$, hence by $B$, which implies that it is a 
Schubert variety $X_t$ since there are only finitely many $B$-orbits. 
The second claim follows immediately: there is only one one-dimensional
Schubert variety $X_s$ (a line), and since the Schubert subvarieties of 
$LG(3,\AA_\CC)$ contained in $\APC$ are exactly the cones over the Schubert
subvarieties of the rational homogeneous variety $\APC$, the interval
$I=[s,t]$ is isomorphic with the Hasse diagram of $\APC$. 
Finally, the third claim was first observed in \cite{cmp2} in connection 
with certain unexpected symmetry properties of quantum cohomology. \qed

\medskip
The picture below shows the Hasse diagram of $LG(3,\OO)$, the {\it 
Freudenthal variety}, which is a homogeneous space of exceptional type, 
with automorphism group $Sp(6,\OO_\CC)$ of type $E_7$. The interval $I$ is in blue 
while $\pi (I)$ is in red. 

\medskip

 \begin{center}
\psset{unit=2.3mm}
\psset{xunit=2.3mm}
\psset{yunit=2.3mm}
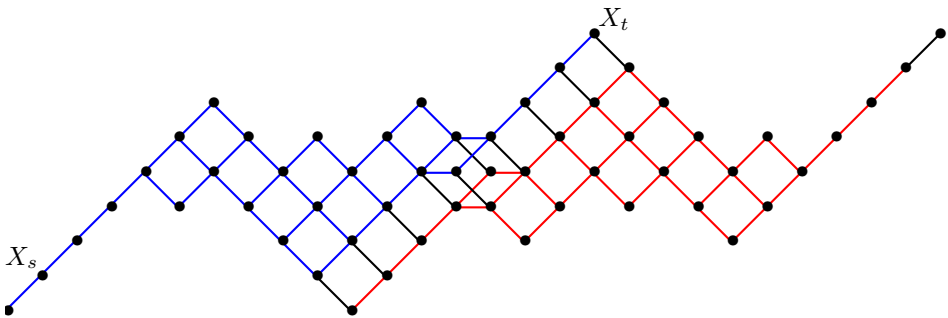
\begin{pspicture*}(50,20)(-8,-4)
\psline(-10,-2)(-8,0)
\psline[linecolor=blue](-8,0)(0,8)
\psline[linecolor=blue](0,8)(4,12)
\psline[linecolor=blue](0,8)(2,6)
\psline[linecolor=blue](2,10)(10,2)\psline(10,2)(12,0)
\psline[linecolor=blue](2,6)(6,10)
\psline[linecolor=blue](6,6)(10,10)
\psline[linecolor=blue](4,12)(12,4)\psline(12,4)(14,2)
\psline[linecolor=blue](10,10)(14,6)\psline(14,6)(16,4)
\psline[linecolor=blue](16,12)(18,10)\psline(18,10)(20,8)
\psline[linecolor=blue](14,10)(16,8)\psline(16,8)(18,6)
\psline[linecolor=blue](8,4)(16,12)
\psline[linecolor=blue](10,2)(18,10)
\psline[linecolor=red](12,0)(20,8)
\psline[linecolor=blue](16,8)(18,8)
\psline[linecolor=red](18,6)(20,6)
\psline[linecolor=blue](18,10)(20,10)
\psline[linecolor=red](20,8)(22,8)
\psline[linecolor=blue](18,8)(26,16)
\psline(18,8)(20,6)\psline[linecolor=red](20,6)(22,4)
\psline[linecolor=red](20,6)(28,14)
\psline[linecolor=red](22,4)(30,12)
\psline(20,10)(22,8)\psline[linecolor=red](22,8)(24,6)
\psline(22,12)(24,10)\psline[linecolor=red](24,10)(28,6)
\psline(24,14)(26,12)\psline[linecolor=red](26,12)(34,4)
\psline(26,16)(28,14)\psline[linecolor=red](28,14)(36,6)
\psline[linecolor=red](28,6)(32,10)
\psline[linecolor=red](32,6)(36,10)
\psline[linecolor=red](34,4)(38,8)
\psline[linecolor=red](36,10)(38,8)
\psline[linecolor=red](38,8)(44,14)\psline(44,14)(46,16)
\multiput(-.2,7.7)(-2,-2){5}{$\bullet$}
\multiput(1.7,9.7)(4,0){5}{$\bullet$}
\multiput(1.7,5.7)(4,0){5}{$\bullet$}
\multiput(3.7,7.7)(4,0){5}{$\bullet$}
\multiput(7.7,3.7)(4,0){3}{$\bullet$}
\multiput(3.7,11.7)(12,0){2}{$\bullet$}
\multiput(9.7,1.7)(4,0){2}{$\bullet$}
\put(11.7,-.3){$\bullet$}
\multiput(17.7,7.7)(4,0){5}{$\bullet$}
\multiput(19.7,9.7)(4,0){5}{$\bullet$}
\multiput(19.7,5.7)(4,0){5}{$\bullet$}
\multiput(21.7,11.7)(4,0){3}{$\bullet$}
\multiput(23.7,13.7)(4,0){2}{$\bullet$}
\multiput(25.7,15.7)(4,0){1}{$\bullet$}
\multiput(21.7,3.7)(12,0){2}{$\bullet$}
\multiput(37.7,7.7)(2,2){5}{$\bullet$}
\put(26.3,16.5){$X_t$}
\put(-8,2.7){$X_s$}
\end{pspicture*}
\end{center}

\medskip\noindent {\it Remark}. There is another connection between
these Hasse diagrams. By Birkhoff's theorem, the Hasse diagram of
$LG(3,\AA_\CC)$, being a distributive lattice, is the lattice of upper
ideals of a poset $P$. This poset is precisely the poset encoded by the
Hasse diagram of $\APC$. In particular the vertices of the latter can
be associated with the join-irreducibles of the former Hasse
diagram. That the Hasse diagram of $LG(3,\AA_\CC)$ is a distributive
lattice is a consequence of the fact that this is a minuscule
homogeneous space \cite{hiller}.

\subsection{Wedge powers of the Jordan algebra}
On the other hand, consider the following problem: decompose the wedge powers
of $\JC$ into irreducible components, with respect to the action of $SL(3,\AA_\CC)$. 
We shall see shortly that this decomposition is multiplicity free. This allows
to define an oriented graph $G(\JC)$ as follows. The vertices are in bijection with the 
components of the wedge powers of $\JC$. There is an edge between a component 
$U$ of $\wedge^k\JC$ and a component $V$ of $\wedge^{k+1}\JC$ if the 
composite map
$$V\otimes\JC^\vee \hookrightarrow \wedge^{k+1}\JC\otimes\JC^\vee
\rightarrow \wedge^{k}\JC\rightarrow U$$
is non-zero. Here the morphism $\wedge^{k+1}\JC\otimes\JC^\vee
\rightarrow \wedge^{k}\JC$ is the natural contraction map, 
and the map $\wedge^{k}\JC\rightarrow U$ is the projection 
with respect to the other irreducible components. 

\begin{prop}
The graph  $G(\JC)$ coincides with the Hasse diagram of $LG(3,\AA_\CC)$. 
\end{prop}

Before giving the proof, we need to recall certain 
properties of the relationship between $\APC$ and $LG(3,\AA_\CC)$. 
First, the latter being minuscule, the Lie algebra $\fg=\fsp(6,\AA)$
of its automorphism group has an associated three-step grading
$$\fg=\fg_{-1}\oplus  \fg_{0}\oplus \fg_{1},$$ 
where $\fg_{0}$ is a reductive Lie algebra with one dimensional center
and with semi-simple part $\fsl(3,\AA_\CC)$, while $\fg_{1}$ is isomorphic
with $\JC$ as a $\fsl(3,\AA_\CC)$-module. The positive part $\fg_{0}\oplus \fg_{1}$
of the grading is the Lie algebra of the parabolic subgroup $P$ of $Sp(6,\AA_\CC)$
such that $G/P=LG(3,\AA_\CC)$. This parabolic is always maximal, hence associated 
to a simple root $\alpha_0$. Its Weyl group $W^P$ can then be defined as 
the stabilizer, inside the Weyl group $W$ of $Sp(6,\AA)$, of the associated 
fundamental coweight $\omega_0^\vee$. The orthogonal hyperplane to $\omega_0^\vee$
cuts the root system $\Phi$ of $Sp(6,\AA)$ along the root system $\Phi_0$
of $SL(3,\AA_\CC)$, which is the root subsystem generated by the simple 
roots except $\alpha_0$. The positive roots which do not belong to $\Phi_0$
are those that appear in $\fg_1$, that is, they are exactly the weights of 
$\JC$. This module is again minuscule, which means that $W^P$ acts transitively
on the roots having positive evaluation on $\omega_0^\vee$.
 
\medskip\noindent {\it Proof of the proposition}. 
We have mentioned the fact that the vertices of the Hasse diagram
are indexed by a set $W_P$, a subset of the Weyl group of $Sp(6,\AA)$. 
This is the set of minimal length representatives of $W/W^P$, and it can
be characterized as follows:
$$W_P=\{w\in W, \quad w(\alpha)\in\Phi^+ \;\forall \alpha\in \Phi_0^+\}.$$
Now, the cotangent space to $LG(3,\AA_\CC)$ at the point $p$ is nothing else
than $\JC$, not only as a vector space but as a $P$-module, hence 
also as an $SL(3,\AA_\CC)$-module since $SL(3,\AA_\CC)$ is  the
semi-simple part of $P$. (In fact the action of the unipotent radical 
of $P$ is trivial, because $LG(3,\AA_\CC)$ is  {\it minuscule}). 
Under such favourable circumstances, the decomposition of the bundle 
of $k$-forms has been obtained by B. Kostant \cite{kos}:  
$$\wedge^k\JC=\bigoplus_{v\in W_P, \ell(v)=k}V_{\rho-v^{-1}(\rho)},$$
where $\rho$ denotes the half-sum of the positive roots in the 
root system of $Sp(6,\AA_\CC)$. In particular, the vertices of $G(\JC)$
are in bijection with $W_P$, hence with the vertices of the Hasse diagram 
of $LG(3,\AA_\CC)$

There remains to check that the edges are the same. In the Hasse diagram,
consider $u$ of length $k$ and $v$ of length $k+1$.  
There is an edge $u\rightarrow v$ if and only if $v=su$ for some
reflection $s$ in the Weyl group $W$. We claim that this is equivalent to the 
condition that $u^{-1}(\rho)-v^{-1}(\rho)$ is  a weight of $\JC$. 
Admitting this, we conclude the proof as follows. 
If $u^{-1}(\rho)-v^{-1}(\rho)$ is not a weight of $\JC$, 
then $V_{\rho-v^{-1}(\rho)}$ cannot be a 
component of $V_{\rho-u^{-1}(\rho)}\otimes \JC$ by \cite[section 131]{zh}. 
If  $u^{-1}(\rho)-v^{-1}(\rho)$ is a weight of $\JC$, then the fact that 
$V_{\rho-v^{-1}(\rho)}$ is a component of $V_{\rho-u^{-1}(\rho)}\otimes \JC$ 
is a special case of the PRV conjecture, proved in \cite{kumar}. 

There remains to prove our claim. First recall that in the minuscule setting
the strong Bruhat order coincides with the weak Bruhat order \cite[Lemma 1.14]{lw}. This means that
the reflexion $s$ must be the simple reflection $s_i$ associated to a 
simple root $\alpha_i$. Since $\ell(v)=\ell(s_iu)=\ell(u)+1$ we must have 
$u(\alpha_i)>0$. We also need that $v$ belongs to $W_P$, which means that 
any  positive root in the subsystem $\Phi_0$ 
must be sent to a positive root.
Since this is already the case for $u$, and since the only positive root
that $s_i$ sends to a negative root is $\alpha_i$, this is equivalent to 
the condition that the positive root $\beta=u^{-1}(\alpha_i)$ does not belong to $\Phi^+_0$.
Since $\Phi_0$ is the set of roots in $\Phi$ orthogonal to $\omega_0^\vee$, 
our condition can be restated as $\omega_0^\vee(\beta)>0$. But this means
that the root $\beta$ appears in $\fg_1$, hence that it is a weight of $\JC$.  
Since $u^{-1}(\rho)-v^{-1}(\rho)=\beta$, our claim follows. 
\qed

\medskip
A nice consequence is that the interval $I=[s,t]$  defines a
submodule of the exterior algebra of $\JC$, namely 
$$L^k= \bigoplus_{\substack{v\in W_P, \ell(v)=k, \\ s\le v\le t}}V_{\rho-v^{-1}(\rho)}.$$
Our refined version of the main observation is the following:

\begin{prop}
$$f_k^a = \dim L^{k+1}.$$
\end{prop}

\medskip
This is straightforward to check case by case. A conceptual 
proof would probably require an interpretation of the 
generalized Dehn-Sommerville equations in representation
theoretic terms. We have no idea of what could be such an
interpretation.

\medskip
As we already mentioned,  although not always irreducible, $L^k$ has 
only a very small
number of irreducible components. More precisely, it contains at most 
three components, as is apparent on the Hasse diagrams of $LG(3,\AA_\CC)$. 
For $k$ small enough $\wedge^k\JC$ is irreducible, hence by symmetry
$L^k$ is also irreducible when $2a+1-k$ is small. 
A natural question to ask is, when the representation is not irreducible,
whether there is any natural way to split the faces into subsets of the 
corresponding dimensions. 

\medskip\noindent {\bf Maximal faces}. In particular, 
the number of faces of maximal dimension is the dimension of the
irreducible module $L^{2a+1}$. This module can be interpreted as follows. 
Recall that each point of $\APC$ defines an $\AA$-line on the dual plane
in $\PP\JC^\vee$. This $\AA$-line is a quadric of dimension $a$, whose
linear span is projective space of dimension $a+1$. Hence an equivariant map 
$$\pi :\APC\rightarrow G(a+2,\JC^\vee )$$
If this map is of degree $d$, in the sense that $\pi^*\mathcal{O}(1)$ is 
equal to the $d$-th power of the hyperplane line bundle on $\APC$, we get 
a non-zero equivariant map $$H^0(G(a+2,\JC^\vee ),\mathcal{O}(1))=
\wedge^{a+2}\JC\stackrel{\pi^*}{\longrightarrow} H^0(\APC,\mathcal{O}(d))=(\JC^\vee)^{(d)},$$ 
the $d$-th Cartan 
power of $\JC^\vee$. Dualizing, we get an inclusion of $\JC ^{(d)}$
inside $\wedge^{a+2}\JC^\vee = \wedge^{2a+1}\JC$. In fact $d=a/2+1$ and 
$$ L^{2a+1}=\JC ^{(a/2+1)}.$$
(This still makes sense for $a=1$ because the hyperplane class of $\RPC$,
the Veronese surface, is divisible by two.)

In general the number of irreducible
components of $L^k$ behaves as follows:
$$\# \mathrm{irred}\; L^{k+1} = \Bigg\{ \begin{array}{l}
     1 \quad \mathrm{if}\; 0\le k\le \frac{a}{2}-1\; \mathrm{or}\; \frac{3a}{2}+1\le k\le 2a, \\                          
     2 \quad \mathrm{if}\; \frac{a}{2}\le k\le a-1\; \mathrm{or}\; a+1\le k\le \frac{a}{2}, \\
     3 \quad \mathrm{if}\; k=a. 
 \end{array} $$
(There is a strange similarity with the homology of $\AP2$.)

\medskip\noindent {\bf Tightness}. 
Note that $\pi (I)=[\pi (s),\pi (t)]$ where the dimension of the Schubert 
variety $X_{\pi (t)}$ is $a+1$, being complementary to the dimension 
of $X_t$, which is $2a+1$ since it is a cone over $\APC$. Since the 
whole Hasse diagram is the disjoint union of $I, \pi (I)$ and the two
extremities, this implies that 
$$L^k=\wedge^k\JC \qquad \mathrm{for}\quad 1\le k\le a+1.$$
This is the algebraic version of tightness. 

\medskip\noindent {\bf Duality}. 
Duality can also be interpreted in the representation theoretic setting. 
Indeed, there exists an exterior automorphism of $SL(3,\AA)$ exchanging 
the representations $\JC$ and its dual $\JC^\vee$. Since 
$$\wedge^k\JC \simeq \wedge^{3a+3-k}\JC^\vee,$$
the graph $G(\JC)$ has an induced symmetry which can be seen to coincide with 
Poincar\'e duality. Since, once again, the disjoint union of $I$ and $\pi (I)$ 
is the whole Hasse diagram minus its two extremities, we must have 
the identity 
$$\wedge^k\JC = L^k\oplus (L^{3a+3-k})^\vee.$$
This is the algebraic version of duality. 
Indeed, taking dimensions, we conclude that
the number of $(k-1)$-dimensional faces is equal to the number of 
$(3a+2-k)$-dimensional ``non-faces''. 
 
\medskip\noindent {\bf Secant defectivity}.
 In \cite{epw}  the minimal triangulation  $\Delta$ of $\RP2$ is considered. 
The Stanley-Reisner ideal $I_\Delta$ defines an arrangement 
of $10$ hyperplanes in $\PP^5$. Over a field $k$ of characteristic
two, the corresponding scheme is Gorenstein
and its canonical bundle is $2$-torsion. Moreover, this scheme
can be flatly deformed into a family of special smooth Enriques 
surfaces in $\PP^5$. 
This family is defined in terms of Lagrangian subspaces in 
$\wedge^3k^6$, endowed with the quadratic form (characteristic
two !) induced by the wedge product. 

In terms of representations (and over $\CC$),
the relevant property is that 
$$\wedge^3(Sym^2\CC^3)=Sym^3\CC^3 \oplus (Sym^2\CC^3)^\vee$$
where both components are Lagrangian (with respect to the 
skew-symmetric form induced by the wedge product). 
In the other cases we have the following substitute:

\begin{prop}
Consider $(L^{2a+1})^\vee$ as an irreducible component of 
$\wedge^{a+2}\JC\simeq\wedge^{2a+1}\JC^\vee$. Then the natural map 
$$(L^{2a+1})^\vee\otimes  (L^{2a+1})^\vee\rightarrow \wedge^{2a+4}\JC=\wedge^{a-1}\JC^\vee$$
is zero. 
\end{prop}

Otherwise stated, $(L^{2a+1})^\vee$ is isotropic with respect to a whole 
system of bilinear forms parametrized by  $\wedge^{a-1}\JC$, which is a very
strong property. Moreover these forms are symmetric for $a\ge 2$, 
exactly as for $a=1$ in characteristic two (but skew-symmetric for 
$a=1$ in characteristic zero...).

\proof Since $(L^{2a+1})^\vee$ is irreducible, it is enough to prove that 
$\omega\wedge\omega'=0$ when $\omega, \omega'$ are two highest weight
vectors. Since $L^{2a+1}$ is a Cartan power of $\JC$, these highest
weight vectors correspond to two points $p, p'$ of the dual $\APC$.
Moreover, we have seen that the associated points in $\wedge^{a+2}\JC$
correspond to the linear spans of the $\AA$-lines on  $\APC$ 
defined by $p$ and $p'$. But two such $\AA$-lines always meet 
non-trivially (we are dealing with a plane projective geometry!), 
and this implies that $\omega\wedge\omega'=0$.\qed

\section{Complexes}

\subsection{A subcomplex of the Koszul complex}
Recall that the wedge powers of $\JC$ can be put together into 
a Koszul complex: for any non-zero linear form $\phi\in\JC^\vee$, 
the contraction by $\phi$, 
$$ \cdots\rightarrow \wedge^{k+1}\JC\stackrel{\phi}{\longrightarrow} 
\wedge^k\JC\rightarrow\cdots $$
defines an exact complex $K^{\bullet}(\phi)$. By their very definition, 
the contraction map by any linear form $\phi$ maps $L^{k+1}$ to $L^k$
and we get a subcomplex $L^{\bullet}(\phi)$
$$ 0\rightarrow L^{2a+1}\rightarrow\cdots\rightarrow L^{k+1}\stackrel{\phi}{\longrightarrow} 
L^k\rightarrow\cdots \rightarrow L^{1}\rightarrow 0.$$
This complex is not exact. Indeed, suppose that $\phi\in\JC^\vee$ is general, 
in the sense that it does not belong to the determinantal hypersurface. 
The stabilizer $SO(\phi)$ of $\phi$ in $SL(3,\AA)$ is then a conjugate of
$SO(3,\AA)=Aut(\JC)$ (such that $\phi$ becomes the identity of the twisted Jordan 
structure). The complex $L^{\bullet}(\phi)$ is $SO(\phi)$-equivariant. 
In particular we consider its Euler characteristic as an element
of the representation ring of $SO(\phi)$. 
A direct check with LiE \cite{LiE} yields:

\begin{prop}\label{inv}
The Euler characteristic of the complex $L^{\bullet}(\phi)$ is
$$\chi(L^{\bullet}(\phi))=\chi_{top}(\AP2)\;[\CC],$$
where $[\CC]$ denotes the class of the trivial representation of $SO(\phi)$.
In particular the Euler characteristic is $SO(\phi)$-invariant. 
\end{prop}

One can also check that 
the $SO(\phi)$-invariants of the complex are
$$(L^{k+1}(\phi))^{SO(3,\AA)}=\Bigg\{ \begin{array}{ll} 
                                  \CC & \mathrm{if}\; k=0,a,2a, \\
                                   0 & \mathrm{otherwise}.
                                 \end{array} $$

The existence of these invariants can be seen as follows. Inside $L^1=\JC$ there
is the invariant hyperplane $\JC_\phi=\phi^\perp$. This is an irreducible $SO(\phi)$-module,
and therefore it admits a unique invariant supplementary line $\ell_\phi\subset \JC$. 
Since $\JC=\JC_\phi\oplus\ell_\phi$ as $SO(\phi)$-modules, we have for any integer $k$  
$$S^k\JC=\oplus_{\ell=0}^k S^\ell\JC_\phi.$$
It turns out that a similar statement holds for Cartan powers:
$$\JC^{(k)}=\oplus_{\ell=0}^k \JC_\phi^{(\ell)}.$$
In particular there is always a unique line of $SO(\phi)$-invariants inside $\JC^{(k)}$, 
hence inside $ L^{2a+1}=\JC ^{(a/2+1)}$. Moreover this line is contained in 
$\wedge^{2a+1}\JC_\phi$, and since $\JC_\phi$ is self-dual of dimension $3a+2$,
there is an induced line of $SO(\phi)$-invariants inside $\wedge^{a+1}\JC_\phi
\subset \wedge^{a+1}\JC=L^{a+1}$.


\medskip
Proposition \ref{inv} suggests the following conjecture:

\medskip\noindent {\bf Conjecture}. {\it 
Let $\phi\in\JC^\vee$ be general. Then 
the inclusion of $L^{\bullet}(\phi)^{SO(\phi)}$ inside 
$L^{\bullet}(\phi)$ is a quasi-isomorphism.}

\medskip
Proposition \ref{inv} also 
shows that $L^{\bullet}(\phi)$ has one of the main properties of 
the face complex of a triangulation of $\AP2$. 

\subsection{The main conjecture}
Let $\Delta$ be a simplicial complex. Associate to each vertex
$v$ of $\Delta$ a variable $x_v$. Let $I_\Delta\subset 
k[x_v, v\in \Delta_0]$  denote the ideal generated by all the square-free monomials 
$x_{v_1}\cdots x_{v_r}$ such that $(v_1,\ldots , v_r)$ is not a face of $\Delta$. 
Then $R=k[x_v, v\in \Delta_0]/I_\Delta$ is the {\it Stanley-Reisner} ring of $\Delta$
\cite{bh}.
When $\Delta$ is a spherical complex, $R$ is a Cohen-Macaulay ring. If $\Delta$ is a 
triangulation of a topological manifold (not necessarily a sphere), then 
$R$ is only Buchsbaum \cite{ns}. 

The {\it face complex} $C_\Delta^\bullet$ is defined by 
$$C_\Delta^k = \bigoplus_{v_1,\ldots , v_k}R_{x_{v_1}\cdots x_{v_k}},$$
where the sum is over all $(k-1)$-dimensional faces. (In order to define 
the morphisms one has to chose an ordering of $\Delta_0$.)
This complex computes the local cohomology of $R$ at the maximal ideal. 
For Buchsbaum modules the local cohomology is closely connected 
with the socle (see \cite{ns}, in particular Corollary 3.5). 

\medskip\noindent {\it Remark}.
Note that a consequence of $\Delta$ not being Cohen-Macaulay is  that the 
$h$-vector is not symmetric. As explained in \cite{ns}, 
the symmetry can be recovered by changing the $h$-vector into a 
$h''$-vector, the modification taking into account the Betti numbers
of the manifold triangulated by $\Delta$. For $\AP2$ we would get 
$$h''_k=h''_{2a-k}=\binom{a+k+1}{k} \qquad \mathrm{for}\quad 0\leq k\le a.$$
These numbers are the dimensions of the graded part of a Gorenstein 
Artinian ring (\cite{ns}, Conjecture 7.3), and by Macaulay's
theorem one can associate a polynomial $F_a$ to this ring. 
What is the significance of $F_a$?

\medskip
We will now define a variant of the face complex. Considering a space $V$ endowed with a basis $e_v$
indexed by vertices of $\Delta$. We can then define
$$L_\Delta^k = \bigoplus_{v_1,\ldots , v_k}\CC e_{v_1}\wedge\cdots\wedge e_{v_k}
\subset \wedge^kV,$$
the sum being again over all $(k-1)$-dimensional faces. 
Since every subset of a face is a face, each contraction map by a linear form 
$\phi\in V^*$ sends $L_\Delta^k$ to $L_\Delta^{k-1}$. 

\medskip\noindent {\bf Conjecture.} {\it There exists a degeneration 
of $L^{\bullet}$ to $L_\Delta^{\bullet}$, for some 
triangulation $\Delta$ of $\APC$ with $3a+3$ vertices.}

\medskip More precisely, such a degeneration should exist inside 
the Koszul complex of $\JC$, which means that we do not need to care
about the morphisms, but only to prove the existence of a 
degeneration $L^k_t$ of each $L^k$ to $L_\Delta^k$ inside the Grassmannian
parametrizing subspaces of $\wedge^k\JC$ of the same dimension. Of course we 
require that for any $\phi\in V^\vee=\JC^\vee$, the contraction map 
by $\phi$ sends $L^{k}_t$ to $L^{k-1}_t$. It would even be natural
to require that for all $k$, 
$$L^{k}_t=Im (L^{2a+1}_t\otimes \wedge^{2a+1-k}\JC^\vee\rightarrow \wedge^{k}\JC) .$$
Then we would only have to degenerate $L^{2a+1}$, subject to the condition 
that these contractions maps have constant rank. 

\begin{prop}
The conjecture is true for $a=1$.
\end{prop}

\proof In this case we only have a three term complex to deal with:
$$ L^3\rightarrow L^2=\wedge^2V\rightarrow L^1=V.$$
Here $V=S^2U$ for $U$ of dimension three. 
In particular there is only $L^3$ to degenerate in the Grassmannian of ten-dimensional
subspaces of $\wedge^3V$, subject to the condition that the contraction map to $\wedge^2V$
is surjective. Since this is an open condition, we can certainly degenerate it to 
the space $L^3_\Delta$ defined by the classical triangulation $\Delta$ of $\RP2$. \qed

\medskip 
It turns out that something rather special happens. Let $u_1, u_2, u_3$ be a basis of $U$, 
and consider the Borel subgroup of $GL(U)$ defined by this basis. Then $L^3=S_{411}U$
is the submodule of $\wedge^3(S^2U)$ with highest weight vector $u_1^2\wedge u_1u_2 \wedge u_1u_3$
with respect to our Borel subgroup. We denote this vector by $(11)(12)(13)$. A basis of $L^3$, 
consisting in eigenvectors of the maximal torus defined by the basis, can 
then be obtained by applying successively the root vectors associated to the 
opposite of the two simple roots of $sl_3$. We get the following diagram. 

$$
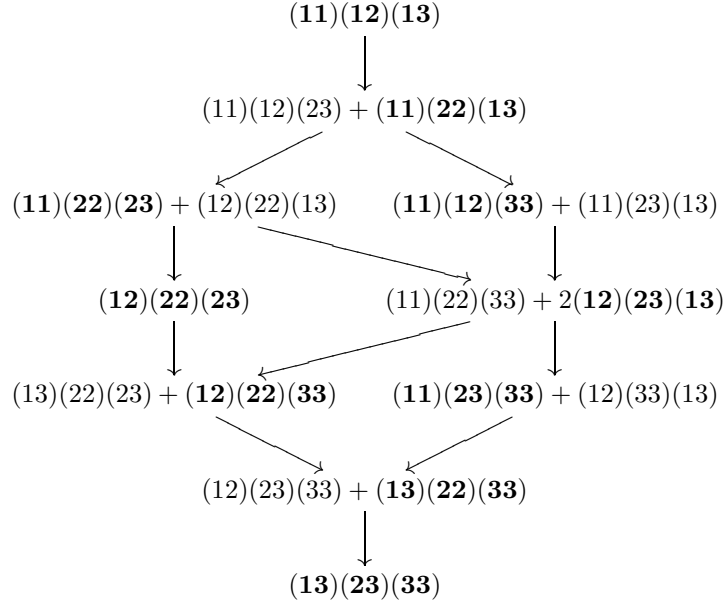
\xymatrix @!0 @C=6pc @R=3pc {
  & {\bf (11)(12)(13)}  \ar[d] \\
  & (11)(12)(23)+{\bf (11)(22)(13)}\ar[rd]\ar[ld] \\
{\bf (11)(22)(23)}+(12)(22)(13)\ar[rrd]\ar[d]  & & {\bf (11)(12)(33)}+(11)(23)(13)\ar[d] \\
 {\bf (12)(22)(23)}\ar[d] & & (11)(22)(33)+2{\bf (12)(23)(13)}\ar[lld]\ar[d] \\
(13)(22)(23)+{\bf (12)(22)(33)}\ar[rd]  &  & {\bf (11)(23)(33)}+(12)(33)(13)\ar[ld]  \\
& (12)(23)(33)+{\bf (13)(22)(33)} \ar[d]  \\
& {\bf (13)(23)(33)}
}$$

\medskip
Now we may consider each basis vector $u_iu_j=(ij)$ of $V$ as a vertex of a triangulation. 
Then a decomposable tensor $(ij)(kl)(mn)$ in $\wedge^3V$ encodes a two-dimensional
face. Not all the vectors in our basis are decomposable, but those that are not
are the sum of only two decomposable vectors, and there is a unique way to choose one 
among these two, for each of the seven non decomposable vectors, in such a way that
the ten faces that we obtain define a triangulation $\Delta$  of $\RP2$
(the minimal triangulation). The terms corresponding to these ten faces 
are indicated in bold on our diagram. In particular we can get a degeneration
of $L^3$ to $L^3_\Delta$ just by rescaling the terms that are not in bold.

\section{The minimal triangulation of $\CP2$, revisited}

In \cite{bk1} the authors exhibited a triangulation of $\CP2$ with nine vertices. 
The list of its 36 maximal faces was obtained with the help of a computer:

$$\begin{array}{lll}
12456\hspace*{10mm} & 45789\hspace*{10mm} & 12378 \\
23456 & 56789 & 12389 \\
13456 & 46789 & 12379 \\
12459 & 34578 & 12678 \\
23567 & 15689 & 23489 \\
13468 & 24679 & 13579 \\
23469 & 35679 & 13689 \\
13457 & 14678 & 12479 \\
12568 & 24589 & 23578 \\
13569 & 34689 & 23679 \\
12467 & 14579 & 13478 \\
23458 & 25678 & 12589
\end{array}$$

The symmetry group $G$ of this triangulation has order 54 and acts 
transitively on the vertices. More specifically, the permutations
$(147)(258)(369)$ and $(123)(456)(789)$ generate a subgroup $H$ 
of the symmetry group isomorphic with $\ZZ_3\times \ZZ_3$, and
this subgroup acts simply transitively on the vertices. Note also
that $H$ has index two in its normalizer $N_G(H)$, which is generated
by $H$ and the involution $\tau=(12)(46)(89)$. The involutions in 
$G$ are all conjugate.  

\medskip
Let us review how this can be connected to our approach. 
For $\AA=\CC$, the Jordan algebra $\JC$ can be identified with the tensor
product $A\otimes B$ of two vector spaces of dimension three. The terms 
of the Koszul complex are given by the Cauchy formula, and the 
subcomplex $L^{\bullet}$ is encoded in the following graph:

$$
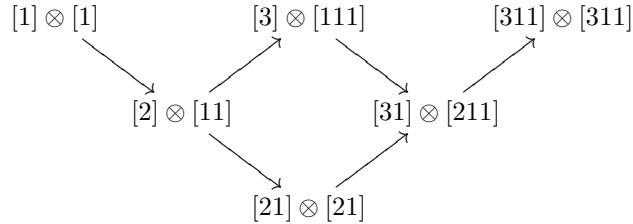
\xymatrix @!0 @C=4pc @R=3pc {
  [1]\otimes [1]  \ar[rd] & & [3]\otimes [111] \ar[rd] & & [311]\otimes [311] \\
 & [2]\otimes [11]  \ar[ru]\ar[dr] & & [31]\otimes [211] \ar[ru] & \\
 & & [21]\otimes [21]\ar[ru] & &
}
$$

\medskip
Our notation here is the following: by $[\mu]\otimes [\nu]$ we mean the tensor
product of Schur powers $S_\mu A\otimes S_\nu B$, plus the symmetric term
$S_\nu A\otimes S_\mu B$ if $\mu\ne\nu$. We have $L^k=\wedge^k(A\otimes B)$
for $1\le k\le 3$, corresponding to the first three columns of the complex. 
On the extreme right, $L^5=S_{311}A\otimes S_{311}B=S^2A\otimes \det A \otimes
S^2B\otimes \det B$ has dimension $36$.  

If we choose a basis $a_1,a_2,a_3$ of $A$ and a basis $b_1,b_2,b_3$ of $B$,
we get a basis $a_i\otimes b_j=(ij)$ of $A\otimes B$ and an induced basis 
of its wedge powers. Note that $L^5=[311]\otimes [311]$ is a multiplicity
free module. As a submodule of $\wedge^5(A\otimes B)$, it is generated
by the highest weight vector $(11)(12)(13)(21)(31)$. Taking into account 
the action of the Weyl group $W=S_3\times S_3$, we get nine decomposable
vectors. Our principle is that each weight vector $(ij)$ should be
identified with a vertex of the triangulation, and each decomposable
vector to a face of this triangulation. Starting from the configuration 
of the nine faces that have to be associated to the nine decomposable
vectors, we are led to the following identification between our weight 
vectors and the vertices of the Brehm-K\"uhnel triangulation:

$$\begin{array}{ccccccccc} 
  1 & 2 & 3 & 4 & 5 & 6 & 7 & 8 & 9 \\
(23) & (32) & (11) & (21) & (33) & (12) & (22) & (31) & (13)
\end{array}$$
 
\medskip
We can then make several observations:
\begin{enumerate}
 \item The transitive action of $H$ on the vertices coincides with the 
natural action of the subgroup $A_3\times A_3$ of the Weyl group $S_3\times S_3$.
 \item The involution $\tau$ coincides with the external 
symmetry $(ij)\mapsto (ji)$. 
 \item The maximal faces split into 
four $H$-orbits of nine elements, corresponding to the four $W$-orbits among
the weights of $L^5$.
 \item Each weight space is generated by a vector which is the sum of one,
two or four decomposable vectors, and exactly one of these decomposable
vectors correspond to a maximal face of the triangulation.
\end{enumerate}
The minimal triangulation of $\CP2$ and its symmetries thus become
much more transparent when interpreted in our representation theoretic setting. 

\section{The quaternionic case, revisited}

In \cite{bk2}, three combinatorial triangulations were constructed of 
a manifold ``like a quaternionic projective plane''. One of these 
triangulations is more symmetric than the two others: its automorphism
group, the icosahedral group $A_5$, acts transitively on the $15$
vertices. It can be characterized as the unique tight
triangulation of a manifold with this symmetry property. 
The authors conjectured that the underlying manifold is really
the quaternionic projective plane $\HP2$, but up to our 
knowledge this conjecture remains open. 

For $\AA=\HH$, the Jordan algebra $\JH$ can be identified with the second
wedge power $\wedge^2A$ of a vector space $a$ of dimension six. The 
subcomplex $L^{\bullet}$ of the Koszul complex is encoded in the following graph:

$$
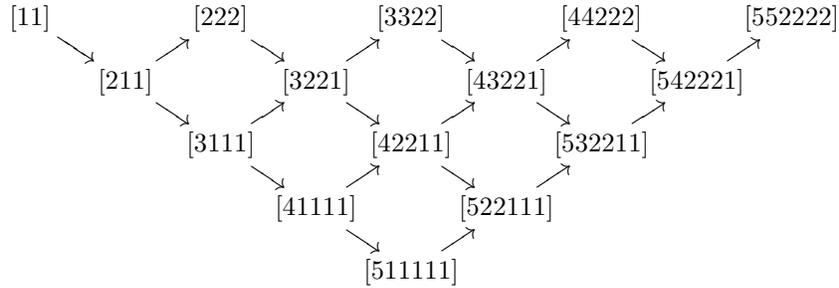
\xymatrix @!0 @C=3pc @R=2pc {
  [11]\ar[rd] & &  [222]\ar[rd]  & & [3322]\ar[rd] & &  [44222] \ar[rd] & &  [552222]\\
 & [211]\ar[rd]\ar[ru] & & [3221]\ar[rd]\ar[ru] & & [43221]\ar[rd]\ar[ru]& & [542221]\ar[ru] & \\
 & & [3111] \ar[rd]\ar[ru] & &  [42211]\ar[rd]\ar[ru] & & [532211]\ar[ru] & &\\
 & & & [41111] \ar[rd]\ar[ru] & & [522111]\ar[ru] & & & \\
 & & & & [511111] \ar[ru] &&&&
}
$$

Here again we denote by $[\lambda]$ the Schur power $S_\lambda A$. In particular $[11]$ 
denotes the minuscule representation $\wedge^2A$, whose fifteen weights (with respect 
to some fixed maximal torus) should represent the fifteen vertices of the Brehm-K\"uhnel triangulation. 
This suggests in particular that the action of the icosahedral group $A_5$ on these vertices, 
which is produced in \cite{bk2} by exhibiting an explicit embedding in $S_{15}$, 
is in fact induced by a much more simple embedding in $S_{6}$. 

This is indeed the case, and we consider this fact as a strong hint that our insights should 
be correct. Consider the permutations defined in cyclic notation by
$$\begin{array}{lll}
p & = & (1)(23456), \\
r_1 &=& (1)(2)(36)(45), \\
s &=& (156)(243), \\
r_2 &=& (4)(5)(36)(12). \\
\end{array}$$
There is an induced action on the set of pairs of distinct integers, 
that we put in correspondence with integers between 1 and 15 by
identifying the following tables:

$$\begin{array}{lllllcl}
45 &36 &12 &\hspace*{2cm} &1  &6  &11  \\
56 &24 &13 &              &2  &7  &12  \\
26 &35 &14 & \qquad \mathrm{vs}             &3  &8  &13  \\
23 &46 &15 &              &4  &9  &14  \\
34 &25 &16 &              &5  &10  &15 
\end{array}$$

It is then straightforward to check that the resulting permutations 
of $S_{15}$ coincide with the permutations denoted $P, R_1, S, R_2$
in \cite{bk2}, pp. 170-171. 

The maximal simplices defining the Brehm-K\"uhnel triangulation
(see \cite{bk2}, Table 2 p. 174) can then be identified with 
sets of nine pairs of integers. In particular, the simplex
denoted $M_1$ corresponds to  $(12)(13)(14)(15)(16)(23)(24)(25)(26)$.
This expression can be seen as defining a highest weight vector
of the representation $[552222]$ inside $\wedge^9[11]$, in complete 
agreement with our expectations. 

\section{Spherical links}

Assuming that the vertex-transitive action by a symmetry group, which
exists for the known triangulations of $\RP2$ and $\CP2$ and for the
conjectural triangulation of $\HP2$, also exists for the hypothetical
case of $\OP2$, one can consider the link of an arbitrary vertex in
one of these triangulations. This link does not depend on the chosen
vertex, up to isomorphism, and defines a triangulation of a sphere of
dimension $2a-1$.

Knowing the number of simplices in the triangulation of $\AP2$, one
can compute the number of simplices in this triangulated sphere by a
double counting argument. First note that in order to count simplices
of the link, one can study the star instead of the link. Consider now
the set of pairs $(v,f)$ where $v$ is a vertex in the triangulation of
$\AP2$ and $f$ is a simplex in the star of $v$. Every $k$-dimensional
simplex of the triangulation of $\AP2$ belongs exactly to the links of
its $k+1$ elements, hence will appear $k+1$ times in the set of pairs
$(v,f)$. One can then count pairs $(v,f)$ such that $f$ is
$k$-dimensional in two different ways.

The results for the spherical triangulations are  listed below 
by increasing dimensions. The spherical link in $\RP2$ is a pentagonal circle. 
According to
\cite[\S 5]{bk2}, the spherical link in $\CP2$ is a non-polytopal
$3$-sphere, called the Br\"uckner-Gr\"unbaum sphere, and the spherical
link in $\HP2$ is a non-polytopal $7$-sphere.

$$\begin{array}{lllll}
& \RP2\; & \CP2\;\; & \HP2\;\;\; & \OP2 \\
\# \mathrm{vertices} & 5 & 8  & 14 &  26 \\
\# \mathrm{1-dim} & \mathbf{5} & 28 & 91 & 325 \\
 &  & 40 & 364 & 2600 \\
\# \mathrm{3-dim}&   & \mathbf{20} & 1001 & 14950 \\ 
&   &  & 1806 & 65780 \\
&   &    & 1974 & 230230 \\
&   &    & 1176 & 657800 \\
\# \mathrm{7-dim}&   &    & \mathbf{294} & 1562275 \\
&   &    &   & 3087370 \\
&  &&& 4964102\\
&  &&& 6255184\\
&  &&& 5922800\\
&  &&& 4022200\\
&  &&& 1838720\\
&  &&& 505648\\
\# \mathrm{15-dim}&  &&& \mathbf{63206} 
\end{array}
$$

\medskip\noindent {\bf Observation.} {\it The number of maximal faces
  in the spherical triangulation is 
  \begin{equation*}
    \frac{3a+2}{a+2}\binom{2a+1}{a+1}.
  \end{equation*}
  This is the dimension of the irreducible
  representation of the Lie algebra $\fso(a+4)$, whose highest weight
  is $a \omega$, where $\omega$ is the fundamental weight defining the vector 
  representation of dimension $a+4$.}

\medskip
The meaning of this observation remains unclear to us. Moreover we could not 
find similar interpretations for the other numbers of faces. 

\smallskip As a curiosity, one can note that the numbers of maximal faces
also appear in the sequence A129869 of the On-Line Encyclopedia of
Integer Sequences ({\texttt oeis.org}), which counts tilting modules
for quivers of type $D$.

\section{Higher ranks}

Much of what we have explained in the previous section remains true for 
higher rank, that is, for the projective spaces $\APn=\RPn, \CPn, \HPn$
with $n\ge 3$. Their complex versions $\APnC$
are homogeneous under the action of a group $SL(n+1,\AA)$. Moreover 
there exists a bigger homogeneous variety $LG(n+1,\AA)$, with 
automorphism group $Sp(2n+2,\AA)$, such that $\APnC$ can be
identified with the space of lines in $LG(n+1,\AA)$ passing 
through a prescribed point $p$. 

$$\begin{array}{cccc}
   \AA &  \RR & \CC &  \HH \\
\APnC & v_2(\PP^n) & \PP^n\times\PP^n & G(2,2n+2) \\
SL(n+1,\AA) & SL(n+1) & SL(n+1)\times SL(n+1)   & SL(2n+2) \\
 LG(n+1,\AA) & LG(n+1,2n+2) & G(n+1,2n+2) & \mathbb{S}_{2n+2} \\
Sp(2n+2,\AA) & Sp(2n+2) & SL(2n+2) & SO(4n+4)
  \end{array}$$

Exactly as for $n=2$, this allows to embed the Hasse diagram of 
$\APnC$ into that of $LG(n+1,\AA)$, as an interval $I=[s,t]$. 
Moreover the Hasse diagram of  $LG(n+1,\AA)$ is canonically 
identified to the graph constructed from the wedge powers of $\JnC$,
the complexification of the space of Hermitian matrices of size $n+1$ 
with coefficients in $\AA$. These wedge powers are given by 
the following classical formulas. 

\medskip\noindent {\bf $\AA = \RR$}. Then $J_{n+1}(\mathbb{R}_{\mathbb{C}})=Sym^2 U$ 
where $U$  has dimension $n+1$. We have
$$\wedge^k J_{n+1}(\mathbb{R}_{\mathbb{C}}) = \bigoplus_{\substack{|\lambda|=k, \\ h(\lambda)\le n}}S_{d_+(\lambda)}U,$$
where the sum is over {\it strict partitions} $\lambda = (\lambda_1>\cdots 
>\lambda_\ell>0)$ of size $k$ and of height $h(\lambda)=\lambda_1$ at most $n$. 
Moreover $d_+(\lambda)$ is the partition of size $2k$ obtained by 
putting together $\lambda$ and its conjugate:
$$d_+(\lambda)=(\lambda_1, \lambda_2+1,\ldots , \lambda_\ell+\ell-1, \ell^{\lambda_\ell},
 (\ell-1)^{\lambda_{\ell-1}-\lambda_\ell-1},\ldots, 1^{\lambda_1-\lambda_2-1}),$$
where powers mean repetitions.

\medskip\noindent {\bf $\AA = \CC$}. Then $J_{n+1}(\mathbb{C}_{\mathbb{C}})
=U\otimes V$ where $U$ and $V$  have dimension $n+1$. We have
$$\wedge^k J_{n+1}(\mathbb{C}_{\mathbb{C}}) = 
\bigoplus_{\substack{|\lambda|=k, \\ \ell(\lambda), h(\lambda)\le n+1}}
S_{\lambda}U\otimes S_{\lambda^\vee}V,$$
where the sum is over partitions $\lambda = (\lambda_1\ge\cdots 
\ge\lambda_\ell>0)$ of size $k$ and of height $h(\lambda)=\lambda_1$
and length $\ell(\lambda)=\ell$ at most $n+1$. 

\medskip\noindent {\bf $\AA = \HH$}. Then $J_{n+1}(\mathbb{H}_{\mathbb{C}})=\wedge^2 U$ 
where $U$ has dimension $2n$. We have
$$\wedge^k J_{n+1}(\mathbb{H}_{\mathbb{C}}) = \bigoplus_{\substack{|\lambda|=k, \\ h(\lambda)\le 2n}}
S_{d_-(\lambda)}U,$$
where the sum is over {\it strict partitions} $\lambda = (\lambda_1>\cdots 
>\lambda_\ell>0)$ of size $k$ and of height $h(\lambda)=\lambda_1$ at most $2n$. 
Moreover $d_-(\lambda)$ is the conjugate partition to $d_+(\lambda)$.

\medskip In each case the Hasse diagram of $LG(n+1,\AA)$ coincides with the 
graph of partitions with the partial order defined by the inclusion relation. 
The minimal element $s$ of $I$ corresponds to the partition of size one, 
while $t$ corresponds to the partition $(n), (n+1,1^n), (2n+1,2n)$ respectively. 
This yields
$$L^{an+1}=\JnC^{(d)} \qquad \mathrm{where}\quad d=a\frac{n-1}{2}+1.$$
The other terms are then easy to write down explicitly:

\medskip\noindent {\bf $\AA = \RR$}. Then $L^k = S_{k+1,1^{k-1}}U$ for 
$1\le k\le n+1$. This implies that  
$$f_{k}^{1,n}=\frac{1}{2}\
\binom{n+k+2}{k+1}\binom{n+1}{k+1}.$$ 

\medskip\noindent {\bf $\AA = \CC$}. Here $L^k = \oplus_{i+j=k-1}S_{i+1,1^j}U
\otimes S_{j+1,1^i}U$ for $1\le k\le 2n+1$. Therefore  
$$f_{k}^{2,n}=\Big(\frac{n+1}{k+1}\Big) ^2
\sum_{i+j=k}\binom{n+i+1}{i}\binom{n}{i}\binom{n+j+1}{j}\binom{n}{j}.$$ 

\medskip\noindent {\bf $\AA = \HH$}. Here $L^k = \oplus_{i+j=k, i>j}S_{d_-(i,j)}U$.

\medskip If we fix a generic element $\phi\in \JnC^\vee$, its stabilizer is 
a conjugate of $SO(n+1,\AA)=Aut(\JnC)$. We expect that the complex $L^\bullet(\phi)$
should be quasi-isomorphic with the complex of $SO(n+1,\AA)$-invariants, with trivial
arrows. Moreover, we have:

\begin{prop}
$$(L^{k+1})^{SO(n+1,\AA)}=\Bigg\{ \begin{array}{ll} 
                                  \CC & \mathrm{if}\; k=0,a,\ldots , na, \\
                                   0 & \mathrm{otherwise}.
                                 \end{array} $$
\end{prop}

\proof Consider for example the case where $a=4$. It follows from the
branching rules from $SL$ to $Sp$ \cite{litt} that a Schur module $S_\mu
U$ has a $Sp(2n)$-invariant if and only if the conjugate partition
$\mu^\vee$ has only even parts, in which case this invariant is unique
up to scalars. For $\mu = d_-(i,j)$, hence $\mu^\vee = d_+(i,j)$, this
means that $i=j+1$ and $j$ is even. Hence $i+j-1=2j$ must be divisible
by four, and the claim follows. \qed

\medskip There is therefore an intriguing relation between these modules
and the cohomology of $\APn$, confirmed by the following statement:

\begin{prop}
For any $a=1,2,4$ and any $n\ge 2$, one has
$$\sum_{k=0}^{na}(-1)^{k}f_{k}^{a,n}=\Bigg\{ \begin{array}{ll} 
                                  \frac{1+(-1)^n}{2} & \mathit{if}\; a=1, \\
                                   n+1 & \mathit{if}\; a=2\; \mathrm{or}\; a=4.
                                 \end{array} $$
 \end{prop}

\medskip An optimistic guess would be that some degeneration of the complex 
    $L^\bullet(\phi)$ should be the Stanley-Reisner complex of some
triangulation of $\APn$. This triangulation would have $f_k^{a,n}$
faces of dimension $k$, in particular it would have exactly $a\binom{n+1}{2}+n+1$
vertices. Unfortunately, this is definitely over-optimistic: it was proved
in \cite{am} that $\mathbb{RP}^3$ does not admit any triangulation 
with only 10 vertices!

\end{document}